\documentclass{article}

\usepackage{amsthm}
\usepackage{amssymb}
\usepackage{amsmath}

\newtheorem{theorem}{Theorem}[section]

\newtheorem{lemma}[theorem]{Lemma}
\theoremstyle{remark}

\theoremstyle{definition}

\theoremstyle{definition}

\numberwithin{equation}{section}

\begin{document}

\title{Existence and Regularity of Optimal Solution\\
for a Dead Oil Isotherm Problem\footnote{Research Report
CM06/I-27, University of Aveiro, 2006. Accepted for publication in
the journal \emph{Differential Geometry -- Dynamical Systems}
(22-July-2006).}}

\author{Moulay Rchid Sidi Ammi
\and Delfim F. M. Torres}

\date{}

\maketitle


\begin{abstract}
We study a system of nonlinear partial
differential equations resulting from the traditional modelling of oil
engineering within the framework of the mechanics of a continuous
medium. Existence and regularity of the optimal solutions for this system
is established.
\end{abstract}


\smallskip


\smallskip


\smallskip


\medskip


\section{Introduction}

We are interested to the existence and regularity of optimal
solution for the following ``dead oil isotherm'' problem:
\begin{equation}
\label{P}
\begin{cases}
\partial_t u - \Delta \varphi(u) = div\left(g(u) \nabla p\right)
& \text{ in } Q_T = \Omega \times (0,T) \, , \\
\partial_t p - div\left(d(u) \nabla p\right) = f
& \text{ in } Q_T = \Omega \times (0,T) \, , \\
\left.u\right|_{\partial \Omega} = 0 \, , \quad \left.u\right|_{t=0} = u_0 \, , \\
\left.p\right|_{\partial \Omega} = 0 \, , \quad \left.p\right|_{t=0}
= p_0 \, ,
\end{cases}
\end{equation}
where $\Omega$ is an open bounded domain in $\mathbb{R}^2$ with a
sufficiently smooth boundary.

Equations \eqref{P} serve as a model of an incompressible biphasic
flow in a porous medium, with applications in the industry of
exploitation of hydrocarbons. To understand the optimal control
problem we consider here, some words about the recovery of
hydrocarbons are in order. At the time of the first run of a layer,
the flow of the crude oil towards the surface is due to the energy
stored in the gases under pressure or in the natural hydraulic
system. To mitigate the consecutive decline of production and the
decomposition of the site, water injections are carried out,
before the normal exhaustion of the layer. The water is injected
through wells with high pressure, by pumps specially drilled with
this end. The pumps allow the displacement of the crude oil towards
the wells of production. The wells must be judiciously distributed,
which gives rise to a difficult problem of optimal control: how to
choose the best installation sites of the production wells?
The cost functional to be minimized comprises all the important
parameters that intervene in the processes.

Existence and uniqueness to the system \eqref{P}, for the case
when the term $\partial_t p$ is missing but for more general
boundary conditions, is established in \cite{mad}.
Optimal control of systems governed by partial differential
equations is investigated in literature by many authors, we can
refer to \cite{Bodart,LeeShilkin,Lions71}. To study existence and
regularity of solutions which provide G\^{a}teaux differentiability
of the nonlinear operator corresponding to \eqref{P}, we are forced to
assume more regularity on the control $f$ as well as to impose
compatibility conditions between initial and boundary conditions.
The considered cost functional comprises four terms and has the form
\begin{equation}
\label{eq:cf} J(u,p,f) = \frac{1}{2} \left\|u - U\right\|_{2,Q_T}^2
+ \frac{1}{2} \left\|p - P\right\|_{2,Q_T}^2 + \frac{\beta_1}{2}
\left\|f\right\|_{2 q_0,Q_T}^{2 q_0} + \frac{\beta_2}{2}
\left\|\partial_t f\right\|_{2,Q_T}^{2}
\end{equation}
where $1 < q_0 < 2$, $\beta_1 > 0$ and $\beta_2 > 0$ are two
coefficients of penalization; $U$ and $P$ are given data. Here $u$
is the reduced saturation of the phase  oil at the moment $t$. The
initial saturation is known and $p$ is the total pressure. The first
two terms in \eqref{eq:cf} make possible to minimize the difference
between the reduced saturation of oil and a given $U$,
respectively the global pressure and a known initial pressure $P$.
We remark that the choice of the objective
functional is not unique. We can always add further terms of
penalization to take into account other properties which one may
want to control. The paper is organized as follows. In Section~\ref{sec2} we
set up the notation, the functional spaces and some important lemmas
used throughout the work. Section~\ref{sec3} is devoted to the existence of optimal
solutions. We obtain necessary estimates on the sequence minimizing
the cost functional which allows us to pass to the limit. Finally,
in Section~\ref{sec4} we establish a regularity theorem.


\section{Notation and Functional Spaces}\label{sec2}

In the sequel we suppose that $\varphi$, $g$ and $d$ are real valued
$C^1$-functions satisfying:
\begin{description}

\item[(H1)] $0 < c_1 \le d(r)$, $\varphi(r) \le c_2$;
$|d'(r)|,\, |\varphi'(r)|,\, |\varphi''(r)|
\leq c_{3}\,  \quad \forall r \in \mathbb{R}$.

\item[(H2)] $u_0$, $p_0$ $\in C^2\left(\bar{\Omega}\right)$, $U$,
$P$ $\in L^2(Q_T)$, where $u_0,\, p_0 : \Omega \rightarrow \mathbb{R}$,
$U,\, P : Q_T \rightarrow \mathbb{R}$, and
$\left.u_0\right|_{\partial \Omega} = \left.p_0\right|_{\partial \Omega} = 0$.

\end{description}

We consider the following spaces:
\begin{equation*}
W_p^{1,0}(Q_T) := L^p\left(0,T,W_p^1(\Omega)\right) = \left\{ u \in
L^p(Q_T), \, \nabla u \in L^p(Q_T) \right\} \, ,
\end{equation*}
endowed with the norm $\left\|u\right\|_{W_p^{1,0}(Q_T)} =
\left\|u\right\|_{p,Q_T} + \left\|\nabla u\right\|_{p,Q_T}$;
\begin{equation*}
W_p^{2,1}(Q_T) := \left\{ u \in W_p^{1,0}(Q_T), \, \nabla^2 u \, ,
\partial_t u \in L^p(Q_T) \right\} \, ,
\end{equation*}
with the norm $\left\|u\right\|_{W_p^{2,1}(Q_T)} =
\left\|u\right\|_{W_p^{1,0}(Q_T)} + \left\|\nabla^2
u\right\|_{p,Q_T} + \left\|\partial_t u\right\|_{p,Q_T}$;
$$
V := \left\{ u \in W_2^{1,0}(Q_T), \, \partial_t u
\in L^2\left(0,T,W_2^{-1}(\Omega)\right) \right\} \, .
$$

We now state some important lemmas that are used later.
Lemma~\ref{lemma2.1} is needed in the proof of our existence result.

\begin{lemma}[\cite{blp}]
\label{lemma2.1} Assume $\Omega \subset \mathbb{R}^n$ is a bounded
domain with a $C^1$-boundary, and a matrix $A(x,t) = \left(A_{i
j}(x,t)\right)$ satisfying the conditions
\begin{equation}
\label{eq:2.1}
\begin{gathered}
\exists \gamma_0 > 0 \text{ such that } A_{i j}(x,t) \xi_i \xi_j
\ge \gamma_0 |\xi|^2 \quad \forall \xi \in \mathbb{R}^n \, , \\
A_{i j} \in L^\infty(Q_T) \, , \quad A_{i j} = A_{j i} \, .
\end{gathered}
\end{equation}
Assume also $f \in L^{2 q_0}(Q_T)$, $u_0 \in W_{2 q_0}^1(\Omega)$
for some $q_0 > 1$ and let $u \in C\left([0,T];\,L^2(\Omega)\right)
\cap W_2^{1,0}(Q_T)$ be a weak solution to the equation
\begin{equation}
\label{eq:2.2}
\begin{gathered}
\partial_t u - div\left(A(x,t) \nabla u\right) = f \text{ in } Q_T \, , \\
\left.u\right|_{\partial \Omega} = 0 \, , \quad \left.u\right|_{t =
0} = u_0 \, .
\end{gathered}
\end{equation}
Then, there exists a constant $q > 1$, depending on $n$, $q_0$,
$\gamma_0$, $Q_T$, and $\left\|A\right\|_{\infty, Q_T}$, such that
$u \in W_{2 q}^{1,0}(Q_T)$ and the estimate
\begin{equation*}
\left\|\nabla u\right\|_{2 q, Q_T} \le c \left( \left\|f\right\|_{2
q, Q_T} + \left\|u_0\right\|_{W_{2q}^1(\Omega)} \right)
\end{equation*}
holds.
\end{lemma}

We use the following two lemmas to get some regularity of weak solutions.

\begin{lemma}[De Giorgi–-Nash–-Ladyzhenskaya–-Uraltseva theorem \cite{lsu}]
\label{lemma2.2} Assume $Q_T = \Omega \times (0,T)$, $\Omega \subset
\mathbb{R}^n$ a $C^1$-bounded domain, and let $f \in L^{s,r}(Q_T) =
L^{s}(0,T,L^r(\Omega))$, $u_0 \in C^\alpha(\bar{\Omega})$ for some
$\alpha_0 > 0$, $\left.u_0\right|_{\partial \Omega} = 0$ and
\begin{equation*}
\label{eq:2.3}
\frac{1}{r} + \frac{n}{2 s} < 1 \, .
\end{equation*}
Assume \eqref{eq:2.1} holds and let $u \in W_2^{1,0}(Q_T)$ be a weak
solution of \eqref{eq:2.2}. Then, there exists $\alpha > 0$ such
that $u \in C^{\alpha,\frac{\alpha}{2}}(\bar{Q}_T)$ and
\begin{equation*}
\label{eq:2.4}
\left\|u\right\|_{C^{\alpha,\frac{\alpha}{2}}(\bar{Q}_T)} \le c
\left(\left\|f\right\|_{L^{s,r}(Q_T)} +
\left\|u_0\right\|_{C^\alpha(\bar{\Omega})}\right) \, .
\end{equation*}
\end{lemma}

\begin{lemma}[\cite{lsu}]
\label{lemma2.3} For any function $u \in
C^{\alpha,\frac{\alpha}{2}}(\bar{Q}_T) \cap L^2\left(0,T; \,
\stackrel{\circ}{W_2^1}(\Omega) \cap W_2^2(\Omega)\right)$ there
exist numbers $N_0$, $\varrho_0$ such that for any $\varrho \le
\varrho_0$ there is a finite covering of $\Omega$ by sets of the
type $\Omega_\varrho(x_i)$, $x_i \in \bar{\Omega}$, such that the
total number of intersections of different $\Omega_{2 \varrho}(x_i)
= \Omega \cap B_{2 \varrho}(x_i)$ does not increase $N_0$. Hence, we
have the estimate
\begin{equation*}
\left\|\nabla u\right\|_{4,Q_T}^4 \le c
\left\|u\right\|_{C^{\alpha,\frac{\alpha}{2}}(\bar{Q}_T)}^2
\mathrm{\varrho}^{2 \alpha} \left( \left\|\nabla^2
u\right\|_{2,Q_T}^2 + \frac{1}{\varrho^2} \left\|\nabla
u\right\|_{2,Q_T}^2\right) \, .
\end{equation*}
\end{lemma}


\section{Existence of Optimal Solution}\label{sec3}

We denote by $(P)$ the problem of minimizing \eqref{eq:cf}
subject to \eqref{P} in the class $\left(u,p,f\right)
\in W_{1}^{2, 1}(Q_{T}) \times V \times L^{2}(Q_{T})$.

\begin{theorem}
\label{theorem3.1} Under hypotheses (H1)-(H2) there is a $q > 1$,
depending on the data of the problem, such that there exists an
optimal solution $\left(\bar{u},\bar{p},\bar{f}\right)$ of
problem $(P)$ verifying:
\begin{gather*}
\bar{u} \in W_{q}^{2,1}(Q_T) \, ,\\
\bar{p} \in C\left([0,T]; L^2(\Omega)\right) \cap W_{2 q}^{1,
0}(Q_T) \, , \quad
\partial_t \bar{p} \in L^2\left(0, T, W_2^{-1}(\Omega)\right) \, , \\
\bar{f} \in L^{2 q_0}(Q_T) \, , \quad \partial_t \bar{f} \in
L^2(Q_T) \, .
\end{gather*}
\end{theorem}

\begin{proof}
Let $(u^{m}, p^{m}, f^{m})\in W_{1}^{2, 1}(Q_{T}) \times V \times
L^{2q_{0}}(Q_{T})$ be a sequence minimizing $J(u, p, f)$. Then we
have
$$  (f^{m}) \mbox{ is bounded in }L^{2q_{0}}(Q_{T}) , $$
$$  (\partial_{t}f^{m}) \mbox{ is bounded in }L^{2}(Q_{T}) . $$
Using the parabolic equation governed by the global pressure $p$ and
Lemma~\ref{lemma2.1}, we know that there exists a number $q>1$ such that
$$ \|\nabla p^{m}\|_{2q, Q_{T}} \leq \left(  \| f^{m}\|_{2q, Q_{T}}+
  \|u_{0}\|_{W_{2q}^{1}(\Omega)} \right). $$
Multiplying the second equation of \eqref{P} by $p$, using the hypotheses
and Young's inequality, we get
  $$
  \sup_{t} \|p^{m}\|^{2}_{2, \Omega}+ \|\nabla p^{m}\|_{2q,
  Q_{T}}^{2} \leq c \|f^{m}\|^{2}_{2,  Q_{T}}.
  $$
  Furthermore, we have that $\partial_{t}p^{m}$ is bounded in
  $L^{2}{(0, T; W_{2}^{-1}(\Omega))}$. By Aubin's Lemma
  \cite{Lions}, $(p^{m})$ is compact in $L^{2}(Q_{T})$.
  Using now the first equation of \eqref{P} we have
  $$  \partial_{t}u^{m} -\varphi'(u^{m})\triangle u^{m} -
  \varphi''(u^{m}) |\nabla u^{m}|^{2}= div(g(u^{m})\nabla p^{m}).$$
  Hence $$ \|u^{m}\|_{W_{q}^{2, 1}(Q_{T})} \leq c, $$
  where all the constants $c$ are independent of $m$. Using
the Lebesgue theorem and the compacity arguments of J. L. Lions
\cite{Lions} we can  extract subsequences, still denoted by $(p^{m}),
(u^{m})$ and $(f^{m})$, such that
$$ p^{m} \rightarrow \overline{p} \mbox{ weakly in } W_{2q}^{1, 0}(Q_{T}),$$
$$ \partial_{t}p^{m} \rightarrow \partial_{t} \overline{p} \mbox{ weakly in }
L^{2}(0, T; W_{2}^{-1}(\Omega),$$
$$ p^{m} \rightarrow \overline{p} \mbox{ strongly in }
L^{2}(Q_{T}),$$
$$ p^{m} \rightarrow \overline{p} \mbox{ a.e. in } L^{2}(Q_{T}),$$
$$ u^{m} \rightarrow \overline{u} \mbox{ a.e. in } L^{2}(Q_{T}),$$
$$ f^{m} \rightarrow \overline{f} \mbox{ weakly in } L^{2q_{0}}(Q_{T}),$$
$$ \partial_{t}f^{m} \rightarrow \partial_{t} \overline{f} \mbox{ weakly in } L^{2}(Q_{T}).$$
The existence of an optimal solution $(\overline{u}, \overline{p},
\overline{f})$ follows, in a standard way, by passing to the limit
in problem \eqref{P} and by using the fact that $J$ is lower
semicontinuous with respect to the weak convergence.
\end{proof}


\section{Regularity of Solutions}\label{sec4}

We now prove some regularity to the solutions
predicted by Theorem~\ref{theorem3.1}.

\begin{theorem}
\label{theorem4.1}
Suppose that (H1) and (H2) are satisfied and
let $\left(\bar{u},\bar{p},\bar{f}\right)$ be an
optimal solution of our problem $(P)$.
Then, there exist a $\alpha > 0$ such that the
following regularity conditions are verified:
\begin{gather}
\bar{p} \in C^{\alpha,\frac{\alpha}{2}}\left(\bar{Q}_T\right) \, , \label{eq:4.1} \\
\bar{u}, \, \bar{p} \in W_{4}^{1,0}(Q_T) \, , \label{eq:4.2} \\
\bar{u}, \, \bar{p} \in W_{2}^{2,1}(Q_T) \, , \label{eq:4.3} \\
\partial_t \bar{u}, \, \partial_t \bar{p}
\in L^\infty\left(0,T; L^2(\Omega)\right) \cap W_2^{1,0}(Q_T) \, , \label{eq:4.4} \\
\bar{u} \in C^{\frac{1}{4}}(Q_T) \, , \label{eq:4.5} \\
\bar{u} \in W_{2 q_0}^{2,1}(Q_T) \, , \quad \bar{p} \in W_{2
q_0}^{2,1}(Q_T) \, . \label{eq:4.6}
\end{gather}
\end{theorem}

\begin{proof}
First, we remark that \eqref{eq:4.1} is an immediate consequence of
Lemma~\ref{lemma2.2}. To show the other results, we begin by proving
the following lemma.
\begin{lemma}
\label{lemma4.2} Consider $(u,p,f)$ solution of \eqref{P}. Assume
that hypotheses (H1) and (H2) hold. Then,
\begin{equation*}
\sup_{t \in (0,T)} \left\|\nabla p\right\|_{2,\Omega}^2 +
\left\|\nabla^2 p\right\|_{2,Q_T}^2 \le c \left( \left\|\nabla
p\right\|_{4,Q_T}^4 + \left\|\nabla u\right\|_{4,Q_T}^4 \right) + c
\end{equation*}
where $c$ depend on $u_0$ and $f$.
\end{lemma}

\begin{proof}
From the second equation of \eqref{P} we have
\begin{equation*}
\partial_t p - d(u) \Delta p = d'(u) \nabla u. \nabla p + f \, .
\end{equation*}
Multiplying this equation by $\partial_t p$ and integrating over
$\Omega$, we obtain
\begin{equation*}
\left\|\partial_t p\right\|_2^2 + \frac{c}{2}
\frac{\partial}{\partial t} \left\|\nabla p\right\|^2 \le c
\int_{\Omega}| \nabla p \nabla u \partial_t p| \, dx
  + \int_{\Omega} |f \partial_t p| \, \, dx \, .
\end{equation*}
Using Young's inequality and integrating in time, we get the desired
estimate.
\end{proof}

To continue the proof of Theorem~\ref{theorem4.1} we need to
estimate $\|\nabla u\|_{4, Q_{T}}$ in function of $\|\nabla p\|_{4,
Q_{T}}$. Then, taking into account the first equation of \eqref{P}, it
is well known that $u \in W_{4}^{1, \frac{1}{2}}(Q_{T})$ and
\begin{equation}\label{eq:4.7}
\|\nabla u\|_{4, Q_{T}} \leq c \|\nabla p\|_{4, Q_{T}}
\end{equation}
(see \cite{ks}). Using  Lemma~\ref{lemma2.3},
we have that for any $\varrho < \varrho_{0}$
\begin{equation*}
\left\|\nabla p\right\|_{4,Q_T}^4 \le c
\left\|p\right\|_{C^{\alpha,\frac{\alpha}{2}}(\bar{Q}_T)}^2
\mathrm{\varrho}^{2 \alpha} \left\{ \left\|\nabla
p\right\|_{4,Q_T}^4 + \frac{1}{\varrho^2} \left\|\nabla
p\right\|_{2,Q_T}^2\right\}+C_{u_{0}, p_{0}, f_{0}} \, .
\end{equation*}
Calling Lemma~\ref{lemma2.2}, we then get \eqref{eq:4.2} for an
eligible choice of $\varrho$. After using \eqref{eq:4.7} we obtain
that $u \in W_{4}^{1, 0}(Q_{T}).$ On the other hand, we have by the
first equation of \eqref{P}  and \eqref{eq:4.2} that $u \in
W_{2}^{2, 1}(Q_{T})$. Moreover, it follows by
Lemma~\ref{lemma4.2} and the fact that $u \in W_{2}^{2, 1}(Q_{T})$
that $p \in W_{2}^{2, 1}(Q_{T})$.

Now, in order to prove \eqref{eq:4.4}, we differentiate both
equations of \eqref{P} with respect to time:
\begin{multline}\label{eq:4.8}
\partial_{tt}u-div\left( \varphi'(u)\nabla \partial_{t}u\right)
-div\left( \varphi''(u)\nabla \partial_{t}u \nabla u\right) \\
= div\left( g'(u) \partial_{t}u \nabla p\right)+ div\left( g(u)\nabla
\partial_{t}p\right),
\end{multline}
\begin{equation}\label{eq:4.9}
\partial_{tt}p-div\left( d(u)\nabla \partial_{t}p\right)
- div\left( d'(u)\nabla \partial_{t}u \nabla p\right)=
 \partial_{t}f.
\end{equation}
Multiplying \eqref{eq:4.9} by $\partial_{t}p$ and integrating over
$\Omega$ we  get
$$ \frac{\partial}{\partial t}\|\partial_{t}p\|^{2}_{2, \Omega}
+ c \|\partial_{t}\nabla p\|^{2}_{2, \Omega} \leq c_{f}+ c
\|\partial_{t}p\|^{2}_{2, \Omega}+c \int_{\Omega}|\partial_{t}u
\nabla p\nabla \partial_{t}p|\,dx. $$
By Young's inequality we have
\begin{equation*}
\begin{split}
\int_{\Omega}|\partial_{t}u \nabla p\nabla \partial_{t}p|\,dx
& \leq \|\partial_{t}u \nabla p\|_{2, \Omega} \|\partial_{t}\nabla p\|_{2, \Omega}\\
& \leq c\|\partial_{t}u \nabla p\|_{2, \Omega}^{2}
+ \frac{c}{2}\|\partial_{t}\nabla p\|_{2, \Omega}^{2}.
\end{split}
\end{equation*}
On the other hand, by Holder's inequality we obtain
\begin{equation*}
\begin{split}
\|\partial_{t}u \nabla p\|_{2, \Omega}^{2}
&= \int_{\Omega}|\partial_{t}u|^{2}|\nabla p|^{2} \\
&\leq \left(\int_{\Omega}|\partial_{t}u|^{4}\right)^{\frac{1}{2}}
 \left(\int_{\Omega}|\partial_{t}p|^{4}\right)^{\frac{1}{2}}=
 \|\partial_{t}u\|_{4, \Omega}^{2} \|\nabla p\|_{4, \Omega}^{2}.
\end{split}
\end{equation*}
Using the following multiplicative inequality \cite{lad}
$$ \|\partial_{t}u\|_{4, \Omega}^{2} \leq c
 \|\partial_{t}u\|_{2, \Omega}
 \|\partial_{t}\nabla u\|_{2, \Omega}  \, \quad
 \forall u \in W_{2}^{1}(\Omega),$$
we obtain
\begin{equation*}
\begin{split}
\int_{\Omega}|\partial_{t}u \nabla p\nabla \partial_{t}p|\,dx
&\leq c \|\partial_{t}u\|_{4, \Omega}^{2} \|\nabla p\|_{4, \Omega}^{2}
+ \frac{c}{2}\|\partial_{t}\nabla p\|_{2, \Omega}^{2}\\
&\leq c \|\partial_{t}u\|_{2, \Omega}\|\partial_{t}\nabla u\|_{2,
\Omega}\|\nabla p\|_{4, \Omega}^{2}+ \frac{c}{2}\|\partial_{t}\nabla p\|_{2, \Omega}^{2}\\
&\leq c \|\partial_{t}u\|_{2, \Omega}^{2}\|\nabla p\|_{4, \Omega}^{4}
+ c \|\partial_{t}\nabla u\|_{2, \Omega}^{2} +\frac{c}{2}
\|\partial_{t}\nabla p\|_{2, \Omega}^{2}.
\end{split}
\end{equation*}
Then,
\begin{multline}
\label{eq:4.10}
\frac{\partial}{\partial t}\|\partial_{t}p\|_{2, \Omega}^{2} + c
 \|\partial_{t}\nabla p\|_{2, \Omega}^{2} \\
 \leq c_{f}+ c \|\partial_{t}p\|_{2, \Omega}^{2}
 + c\|\partial_{t}\nabla u\|_{2, \Omega}^{2}
 + c \|\partial_{t} u\|_{2, \Omega}^{2} \|\nabla p\|_{4,
 \Omega}^{4}.
\end{multline}
Multiplying \eqref{eq:4.8} by $\partial_{t}u$ and integrating over
$\Omega$, we get
\begin{multline*}
\frac{\partial}{\partial t}\|\partial_{t}u\|_{2, \Omega}^{2} + c
 \|\partial_{t}\nabla u\|_{2, \Omega}^{2} \\
 \leq - \int_{\Omega}\varphi''(u)\partial_{t}u \nabla u \partial_{t}\nabla
 u - \int_{\Omega}g'(u)\partial_{t}u \nabla p \nabla u
- \int_{\Omega}g(u)\partial_{t} \nabla p \nabla u.
\end{multline*}
Similar as before, we have:
$$
\left|\int_{\Omega}\varphi''(u)\partial_{t}u \nabla u \partial_{t}\nabla u \right|
\leq c \|\partial_{t} u\|_{2, \Omega}^{2} \|\nabla u\|_{4, \Omega}^{4}
 +c \|\partial_{t} \nabla u\|_{2, \Omega}^{2}, $$
$$
\left| \int_{\Omega}g'(u)\partial_{t}u \nabla p \nabla u \right|
\leq c \|\partial_{t} u\|_{2, \Omega}^{2} \|\nabla p\|_{4, \Omega}^{4}
 + c \| \nabla u\|_{2, \Omega}^{2},
 $$
 $$
\left|\int_{\Omega}g(u)\partial_{t} \nabla p \nabla u \right|
\leq c \|\partial_{t} \nabla p\|_{2, \Omega}^{2}+ c \| \nabla u\|_{2, \Omega}^{2}.
$$
It follows by using \eqref{eq:4.7} that
\begin{equation}\label{eq:4.11}
\frac{\partial}{\partial t} \|\partial_{t}  u\|_{2, \Omega}^{2}+ c
\|\partial_{t}  \nabla u\|_{2, \Omega}^{2} \leq c
 \|\partial_{t}
u\|_{2, \Omega}^{2} \| \nabla p\|_{4, \Omega}^{4} + c \|\partial_{t}
\nabla p\|_{2, \Omega}^{2}+ c \| \nabla u\|_{2, \Omega}^{2}.
\end{equation}
Calling \eqref{eq:4.10} and \eqref{eq:4.11} together, it yields
\begin{multline*}
\frac{d}{dt}\left \{ \|\partial_{t} u\|_{2, \Omega}^{2}+
\|\partial_{t} p\|_{2, \Omega}^{2} \right \}+ \|\partial t \nabla
u\|_{2, \Omega}^{2}+
\|\partial t \nabla p\|_{2, \Omega}^{2}\\
\leq c \left( 1+ \| \nabla p\|_{4, \Omega}^{4}\right) \left
\{\|\partial_{t}
 u\|_{2, \Omega}^{2}+ \|\partial_{t} p\|_{2, \Omega}^{2}  \right\}+ c_{u_{0},
p_{0}, f}.
\end{multline*}
We thus obtain \eqref{eq:4.4} by applying Gronwall lemma.
On the other hand, we have
$$
\|\partial_{t} u\|_{4, \Omega}^{2} \leq c \|\partial_{t} u\|_{2,
\Omega} \|\partial_{t} \nabla u\|_{2, \Omega},
$$
and from \eqref{eq:4.4} we obtain $\partial_{t} u \in
L^{4}(Q_{T})$. Using then \eqref{eq:4.2} and the fact that
$W_{4}^{1}(Q_{T}) \hookrightarrow
C^{\frac{1}{4}}(\overline{Q_{T}})$, the regularity
estimate \eqref{eq:4.5} follows.
Finally, the right hand side of the first equation of \eqref{P}
belongs to $L^{4}(Q_{T})\hookrightarrow L^{2q_{0}}(Q_{T})$ as
$2q_{0}\leq 4$. Using thus \eqref{eq:4.3} we get $u \in
W_{2q_{0}}^{2, 1}(Q_{T})$. Since $f \in L^{2q_{0}}(Q_{T})$, the same
estimate follows for $p$ from the second equation of the system \eqref{P}
and we conclude with \eqref{eq:4.6}.
\end{proof}


\subsection*{Acknowledgements}

The authors were supported by FCT (\emph{The Portuguese Foundation
for Science and Technology}): Sidi Ammi through the fellowship
SFRH/BPD/20934/2004, Torres through the R\&D unit CEOC of the
University of Aveiro.



\bigskip

\noindent
Moulay Rchid Sidi Ammi and Delfim F. M. Torres\\
Department of Mathematics, University of Aveiro\\
3810-193 Aveiro, Portugal\\
\texttt{\{sidiammi, delfim\}@mat.ua.pt}

\end{document}